\documentclass[12pt]{amsart}
\usepackage{latexsym}

\newtheorem{thm}{Theorem}[section]
\newtheorem{propn}[thm]{Proposition}
\newtheorem{lemma}[thm]{Lemma}
\newtheorem{coro}[thm]{Corollary}
\newtheorem{rem}[thm]{Remark}

\begin{document}

\baselineskip=16pt

\title{Deformations of the Picard bundle}

\author{I. Biswas}

\address{School of Mathematics, Tata Institute of
Fundamental Research, Homi Bhabha Road, Bombay 400005, India}

\email{indranil@math.tifr.res.in}

\author{L. Brambila-Paz}

\address{CIMAT, Apdo. Postal 402, C.P. 36240. Guanajuato, Gto,
M\'exico}

\email{lebp@fractal.cimat.mx}

\author{P. E. Newstead}

\address{Department of Mathematical Sciences, The University of
Liverpool, Peach Street, Liverpool, L69 7ZL, England}

\email{newstead@liverpool.ac.uk}

\subjclass{14H60, 14J60}

\thanks{All authors are members of the research group VBAC
(Vector Bundles on Algebraic Curves), which is partially supported
by EAGER (EC FP5 Contract no. HPRN-CT-2000-00099) and by EDGE (EC
FP5 Contract no. HPRN-CT-2000-00101).  The second author acknowledges
the support of CONACYT grant 28492-E}

\date{}

\begin{abstract}

Let $X$ be a nonsingular algebraic curve of genus $g\geq 3$, and let
${\mathcal M}_{\xi}$ denote the moduli space of stable
vector bundles of rank $n\ge2$ and degree $d$ with fixed determinant
$\xi$ over $X$ such that $n$ and $d$ are
coprime and $d> n(2g-2)$. We assume that if $g=3$ then $n\geq 4$
and if $g=4$ then $n\geq 3$. Let ${\mathcal W}_{\xi}(L)$ denote
the vector bundle over ${\mathcal M}_\xi$ defined by the direct
image $p_{\mathcal M_{\xi} *}({\mathcal U}_{\xi} \otimes p^*_X L)$
where
${\mathcal U}_{\xi}$
is a universal vector bundle over $X\times {\mathcal
M}_{\xi}$ and $L$ is a line bundle over $X$ of degree zero.
The space of infinitesimal deformations of ${\mathcal
W}_{\xi}(L)$ is proved to be isomorphic to $H^1(X,\,{\mathcal O}_X)$.
This construction gives a complete family
of vector bundles over ${\mathcal M}_\xi$
parametrized by the Jacobian $J$ of $X$
such that ${\mathcal W}_{\xi}(L)$ is the
vector bundle corresponding to ${L} \in J$.
The connected component of the moduli space of stable sheaves with the
same Hilbert polynomial as ${\mathcal W}_{\xi}({\mathcal O})$
over ${\mathcal M}_{\xi}$
containing  ${\mathcal W}_{\xi}({\mathcal O})$
is in fact isomorphic to $J$ as a polarised variety.
\end{abstract}

\maketitle

\section{Introduction}

Let $X$ be a connected
nonsingular projective algebraic curve of genus $g \geq 2$ defined
over the complex numbers.
Let
$J$ denote the Jacobian (Picard variety) of $X$ and $J^d$ the variety
of line bundles of degree $d$ over $X$; thus in particular $J^0=J$.
Suppose $d\ge 2g-1$ and let $\mathcal L$ be a Poincar\'e (universal)
bundle over $X\times J^d$.  If we denote by $p_J$ the natural
projection from $X\times J^d$ to $J^d$, the direct image
$p_{J*}\mathcal L$ is then locally free and is called the Picard
bundle of degree $d$.

These bundles have been investigated by a number of authors over at
least the last 40 years. It may be noted that the projective bundle
corresponding to $p_{J*}\mathcal L$ can be identified with the
$d$-fold symmetric product $S^d(X)$. Picard bundles were studied in
this light by  A. Mattuck \cite{ma1, ma2} and I. G. Macdonald
\cite{mc} among others; both Mattuck and Macdonald gave formulae for
their Chern classes. Somewhat later R. C. Gunning \cite{gu1, gu2} gave
a more analytic treatment involving theta-functions. Later still, and
of especial relevance to us, G. Kempf \cite{ke} and S. Mukai
\cite{mu1} independently studied the deformations of the Picard
bundle; the problem then is to obtain an inversion formula showing
that all deformations of
$p_{J*}\mathcal L$ arise in a natural way. Kempf and Mukai proved
that $p_{J*}\mathcal L$ is simple and that, if $X$ is not
hyperelliptic, the space of infinitesimal deformations of
$p_{J*}\mathcal L$ has dimension given by
$$
\dim\, H^1(J^d, \, End
(p_{J*}{\mathcal L}))=
2g.
$$
Moreover, all the infinitesimal deformations arise from genuine
deformations. In fact there is a complete family of deformations of
$p_{J*}\mathcal L$ parametrised by $J\times {\rm Pic}^0(J^d)$, the
two factors corresponding respectively to translations  in $J^d$ and
deformations of $\mathcal L$ (\cite[\S9]{ke}, \cite[Theorem
4.8]{mu1}). (The deformations of $\mathcal L$ are given by ${\mathcal
L}\mapsto{\mathcal L}\otimes p_J^* L$ for $L\in{\rm
Pic}^0(J^d)$.) Since $J$ is a principally polarised abelian variety
and $J^d\cong J$, ${\rm Pic}^0(J^d)$ can be identified with $J$
(strictly speaking ${\rm Pic}^0(J^d)$ is the dual abelian variety,
but the principal polarisation allows the identification).

Mukai's paper \cite{mu1} was set in a more general context involving a
transform which provides an equivalence between the derived category
of the category of ${\mathcal O}_A$-modules over an abelian variety
$A$ and the corresponding derived category on the dual abelian
variety $\hat A$. This technique has come to be known as the
Fourier--Mukai transform and has proved very useful in studying
moduli spaces of sheaves on abelian varieties and on some other
varieties.

Our object in this paper is to generalise the results of Kempf and
Mukai on deformations of Picard bundles to the moduli spaces of
higher rank vector bundles over $X$ with fixed determinant. In
particular we compute the space of infinitesimal deformations of the
Picard bundle in this context and also identify a complete family of
deformations. While the analogy is not precise, this identification
can be seen as a type of Fourier--Mukai transform.

We fix a holomorphic line bundle $\xi$ over $X$ of degree $d$.
Let ${\mathcal M}_{\xi}\, :=\, {\mathcal M}_{\xi}(n,d)$ be the
moduli space of stable vector bundles $E$ over $X$ with ${\rm
rank}(E) =n\ge2$, ${\rm deg}(E) =d$ and $\bigwedge^n E =\xi$.
We assume that $n$ and $d$ are coprime, ensuring
the smoothness and completeness of ${\mathcal M}_{\xi}$, and that
$g\ge3$. We  assume also that if $g=3$ then  $n\geq 4$ and if $g=4$
then
$n\geq 3$. The case $g=2$ together with the three special cases $g=3$
with
$n=2,3$ and
$g=4$ with
$n=2$ are omitted in our main result since the method of proof
does not cover these cases.

It is known that there is a universal vector bundle over
$X\times {\mathcal M}_{\xi}$. Two such universal bundles
differ by tensoring with the pullback of a line bundle on
${\mathcal M}_{\xi}$. However, since ${\rm Pic}({\mathcal
M}_{\xi}) \, =\, {\mathbb Z}$, it is possible to choose canonically
a universal bundle. Let $l$ be the smallest positive number such
that $ld \equiv 1$ mod $n$. There is a unique universal vector
bundle ${\mathcal U}_{\xi}$ over $X\times {\mathcal M}_{\xi}$
such that ${\bigwedge}^n {\mathcal U}_{\xi}\big\vert_{\{x\}\times
{\mathcal M}_{\xi}} \,=\, {\Theta}^{\otimes l}$ \cite{Ra}, where
$x\in X$ and $\Theta$ is the ample generator of ${\rm
Pic}({\mathcal M}_{\xi})$.  Henceforth, by a universal bundle we
will always mean this canonical one. We denote by $p_X$ and
$p_{\mathcal M}$ the natural projections of $X\times{\mathcal
M}_{\xi}$ onto the two factors.

Now suppose that $d>n(2g-2)$.
For any $L\in J$, let
$$
{\mathcal W}_{\xi}(L)\, :=\, p_{{\mathcal M}*} ({\mathcal
U}_{\xi}\otimes p^*_X L)
$$
be the direct image. The assumption on $d$ ensures that
${\mathcal W}_{\xi}(L)$ is a locally free sheaf on
${\mathcal M}_{\xi}$ and all the higher direct images
of ${\mathcal U}_{\xi}\otimes p^*_X L$ vanish.
The rank of ${\mathcal W}_{\xi}(L)\,$ is $d+n(1-g)$ and
$H^i({\mathcal M}_{\xi},\, {\mathcal W}_{\xi}(L))\,\cong\,
H^i(X\times {\mathcal M}_{\xi}, \,{\mathcal U}_{\xi}\otimes p^*_X
L)$. By analogy with the case $n=1$, we shall refer to the bundles
${\mathcal W}_{\xi}(L)$ as {\em Picard bundles}.

Our first main result concerns the infinitesimal deformations of
${\mathcal W}_{\xi}(L)$. We prove

\medskip
{\bf Theorem \ref{2.9}.}
{\em For any line bundle $L\,\in\, J$,
the space of infinitesimal deformations of the vector bundle
${\mathcal W}_{\xi}(L)$, namely $H^1({\mathcal M}_\xi ,\, End
( {\mathcal W}_{\xi}(L)))$, is canonically isomorphic to
$H^1(X,\, {\mathcal O}_X)$. In particular,
$$
\dim H^1({\mathcal M}_\xi ,\, End  ( {\mathcal W}_{\xi}(L)))
\, =\, g\, .
$$}
\medskip

In the special case where $n=2$ and $L={\mathcal O}_X$, this was
proved by V. Balaji and P. A. Vishwanath in \cite{bv} using a
construction of M. Thaddeus \cite{ta}. For all $n$, it is known that
${\mathcal W}_{\xi}({\mathcal O}_X)$ is simple \cite{bef} and indeed
that it is stable (with respect to the unique polarisation of
${\mathcal M}_{\xi}$)
\cite{bbgn}; in fact the proof of stability generalises easily to
show that ${\mathcal W}_{\xi}(L)$ is stable. In this context, note
that Y. Li \cite{li} has proved a stability result for Picard
bundles over the non-fixed determinant moduli space ${\mathcal
M}(n,d)$, but this does not imply the result for ${\mathcal M}_{\xi}$.
As a byproduct of our proof of Theorem
\ref{2.9}, we obtain a new proof that
${\mathcal W}_{\xi}(L)$ is simple (Corollary \ref{2.8}).

We can consider the bundles $\{ {\mathcal
W}_{\xi}(L)\}$ as a family of bundles over ${\mathcal M}_{\xi}$,
parametrized by
$J$. We prove that this is a complete family of deformations, both
globally and in the local sense that
the infinitesimal deformation map
$$
H^1(X,\, {\mathcal O}_X)\, \longrightarrow\, H^1({\mathcal M}_{\xi},\,
End({\mathcal W}_{\xi}(L)))
$$
is an isomorphism for all $L$ (Corollary \ref{5.2}). As we noted,
the bundles ${\mathcal W}_{\xi}(L)$ are stable.
We denote by ${\mathcal M}^0_X({\mathcal W}_{\xi})$
the connected component of the moduli space of stable
sheaves with the same Hilbert polynomial as
${\mathcal W}_{\xi}({\mathcal O})$ on ${\mathcal M}_{\xi}$
containing ${\mathcal W}_{\xi}({\mathcal O})$. We prove

\medskip
{\bf Theorem \ref{5.1}.}
{\em The morphism
$$
\phi\, : \,J\, \longrightarrow\,
{\mathcal M}^0_X({\mathcal W}_{\xi})
$$
given by $\phi(L)={\mathcal W}_{\xi}(L)$ is an isomorphism of
polarised varieties.}
\medskip

A further consequence is the following  Torelli theorem: if
$X$ and $X'$ are smooth projective curves, $\xi$, $\xi'$ are line
bundles of degree $d$ over $X$, $X'$ respectively and
${\mathcal M}^0_X({\mathcal W}_{\xi})
\cong  {\mathcal M}^0_{X'}({\mathcal W}_{\xi'})$ as polarised
varieties, then
$X\,\cong\, X'$ (Corollary \ref{5.4}).

{\it Notation and assumptions.} We work throughout over the
complex numbers and suppose that $X$ is a connected nonsingular
projective algebraic curve of genus $g\ge3$. We suppose that $n\ge2$
and that if
$g=3$ then $n\ge4$ and if $g=4$ then
$n\ge3$. We assume moreover that $(n,d)=1$ and $d>n(2g-2)$. In
general, we denote the natural projections of
$X\times Y$ onto its factors by $p_X$, $p_Y$, writing $p_{\mathcal
M}$ for
$p_{\mathcal M_{\xi}}$ for simplicity. For a variety $X\times Y\times
Z$, we denote by $p_i\ (i=1,2,3)$ the projection onto the $i$-th
factor and by $p_{ij}$ the projection onto the Cartesian
product of the $i$-th and the
$j$-th factors. Finally, for any $x\in X$, we denote by ${\mathcal
U}_x$ the bundle over ${\mathcal M}_{\xi}$ obtained by restricting
${\mathcal U}_{\xi}$ to $\{x\}\times{\mathcal M}_{\xi}$.

\smallskip
{\small {\it Acknowledgements.} Part of this work was done
during a visit of the authors to ICTP, Italy. The
authors thank ICTP for its hospitality. The third author would like
also to thank the Isaac Newton Institute, Cambridge and the organisers
of the HDG programme for their hospitality during the completion of
work on this paper.}
\smallskip

\section{Cohomology of ${\mathcal W}_{\xi}(L_1)\otimes{\mathcal
W}_\xi(L_2)^*$}

To compute the cohomology groups $H^i({\mathcal M}_{\xi},
{\mathcal W}_{\xi}(L_1)\otimes{\mathcal W}_\xi(L_2)^*)$ we need
the following propositions.

\begin{propn}\label{2.1}
If $L_1$ and $L_2$ are two line bundles in $J$
then
$$
H^i({\mathcal M}_\xi,
\, {\mathcal W}_{\xi}(L_1)\otimes {\mathcal W}_{\xi}(L_2)^* )\, \cong \,
H^i(X\times {\mathcal M}_\xi,\,
{\mathcal U}_{\xi}\otimes p_X^*L_1\otimes p_{\mathcal
M}^*{\mathcal W}_\xi (L_2)^*)
$$
for any $i \geq 0$. Moreover,
$$
H^i({\mathcal M}_\xi ,\,  {\mathcal W}_{\xi}(L_1)\otimes
{\mathcal W}_{\xi}(L_2)^*)
\,\cong \, H^{i+1}(X\times {\mathcal M}_\xi \times X, \, p_{12}^*
{\mathcal
U}_{\xi}\otimes p^*_1L_1 \otimes p^*_{23}{\mathcal U}^*_\xi\otimes
p^*_3 L_2^*\otimes  p^*_3K_X)
$$
for $i \geq 0$, where $K_X$ is the canonical line bundle over $X$.
\end{propn}

\begin{proof}
By the assumption that $d >n(2g-2)$ we have
that $H^1(X,\, {\mathcal
U}_{\xi}\big\vert_{X\times \{v\}}\otimes L_1) =0$ for all $v\in
{\mathcal M}_\xi$.  Using the projection formula and the Leray
spectral sequence we have
$$
H^i(X\times {\mathcal M}_\xi,\,
{\mathcal U}_{\xi}\otimes p^*_X L_1\otimes p_{\mathcal
M}^*{\mathcal W}_\xi (L_2)^*)\, \cong\, H^i({\mathcal M}_\xi,
\, {\mathcal W}_{\xi}(L_1)\otimes {\mathcal W}_{\xi}(L_2)^* )\, .
$$
This proves the first part.

For every stable vector bundle $E$ of rank $n$
and degree $d$,
$$
H^0(X,\, E^*\otimes L^*_2\otimes K_X) \cong
H^1(X,\, E\otimes L_2)^* \, =\, 0\, .
$$
Consequently, the projection formula gives
$$
{\mathcal R}^{i} p_{12*}(p_{23}^*({\mathcal U}^*_\xi)\otimes
p_3^*L^*_2 \otimes p_3^* K_X) \, =\, 0
$$
for $i\not=1$, and we have by relative Serre duality
$$
{\mathcal R} ^{1} p_{12*}(p_{23}^*({\mathcal U}^*_\xi)\otimes
p_3^*L^*_2 \otimes p_3^*K_X)\, \cong\, p_{\mathcal M}^* {\mathcal
W}_\xi (L_2)^*\, ,
$$
Finally, using the projection formula and the Leray spectral
sequence, it follows that
$$
H^i(X\times {\mathcal M}_\xi,\,
{\mathcal U}_{\xi}\otimes p^*_XL_1 \otimes
p_{\mathcal M}^*{\mathcal W}_\xi (L_2)^*).
$$
$$
\cong\,
H^{i+1}(X\times {\mathcal M}_\xi \times X, \, p_{12}^*
{\mathcal U}_{\xi} \otimes p^*_1L_1\otimes
  p^*_{23}{\mathcal U}^*_\xi \otimes p^*_3L^*_2 \otimes p^*_3 K_X)
$$
for $i\geq 0$. Thus,
$$
H^{i+1}(X\times {\mathcal M}_\xi \times X, \, p^*_{12}
{\mathcal U}_{\xi}\otimes p^*_1L_1 \otimes
p^*_{23}{\mathcal U}^*_\xi\otimes p^*_3 L^*_2\otimes p^*_3K_X)\,
\cong\, H^i({\mathcal M}_\xi ,\,   {\mathcal W}_{\xi}(L_1)
\otimes {\mathcal W}_{\xi}(L_2)^*)
$$
as asserted in the proposition.
\end{proof}

\begin{rem}\label{2.2}
\begin{em}
Proposition \ref{2.1} can be formulated
in a more general context.
Let $V_1$, $V_2$ be flat families of vector
bundles over $X$
parametrised by a complete irreducible variety $Y$
such that for each $y\in Y$ we have
$H^1(X,\, V_i|_{X\times \{y\}})\,=\,0$ for $i=1,2$. Under this
assumption
$$
H^i(Y,\, p_{Y*}V_1\otimes (p_{Y*}V_2)^*)\,\cong\, H^{i+1}(X\times Y
\times X,
\, p^*_{12}V_1 \otimes p^*_{23}V_2^*\otimes p^*_3K_X)\, .
$$
The proof is the same as for Proposition \ref{2.1}.
\end{em}\end{rem}

Denote by ${\mathcal R}^i$ the $i$-th direct image of $ p^*_{12}
{\mathcal U}_{\xi}\otimes p^*_1L_1 \otimes
p^*_{23}{\mathcal U}^*_\xi$ for the projection $p_3$, that is
$${\mathcal R}^i := {\mathcal R}^i p_{3*}( p^*_{12}
{\mathcal U}_{\xi}\otimes p^*_1L_1 \otimes
p^*_{23}{\mathcal U}^*_\xi).$$
{}From the Leray spectral sequence we obtain the  exact
sequences
$$
0\,\longrightarrow\, H^1(X,\,{\mathcal R}^0\otimes L^*_2\otimes K_X)
\,\longrightarrow\,
H^{1}(X\times {\mathcal M}_\xi \times X, \, p^*_{12}
{\mathcal U}_{\xi}\otimes p^*_1L_1 \otimes
p^*_{23}{\mathcal U}^*_\xi\otimes p^*_3 L^*_2
\otimes p^*_3K_X)
$$
\begin{equation}\label{1}
\hspace{5.5cm}\, \longrightarrow\, H^0(X,{\mathcal R}^1\otimes L^*_2\otimes
K_X)\,\longrightarrow \,0
\end{equation}
and
$$
0\, \longrightarrow\, H^1(X,{\mathcal R}^1\otimes L^*_2\otimes K_X)
\, \longrightarrow\,
H^{2}(X\times {\mathcal M}_\xi \times X, \, p^*_{12}
{\mathcal U}_{\xi}\otimes p^*_1L_1 \otimes
p^*_{23}{\mathcal U}^*_\xi\otimes p^*_3 L^*_2\otimes p^*_3K_X)
$$
\begin{equation}\label{2}
\hspace{5.5cm}\, \longrightarrow\, H^0(X,{\mathcal R}^2\otimes L^*_2
\otimes K_X)\, \longrightarrow\, 0\, .
\end{equation}

\begin{propn}\label{2.3}
${\mathcal R}^0 \,=\, 0$.
\end{propn}

\begin{proof}
Note that for any $x \in X$
\begin{equation}\label{3}
H^0(X\times {\mathcal M}_{\xi},\, {\mathcal U}_{\xi} \otimes
p^*_XL_1
\otimes p_{\mathcal M}^*{\mathcal U}^*_{x}) \cong H^0(X ,\,
L_1\otimes  p_{X*}({\mathcal U}_{\xi}
\otimes p_{\mathcal M}^*{\mathcal U}^*_{x}))\, .
\end{equation}

{}From \cite{ib} and \cite{nr} we know that for generic $y \in X$,
the two vector bundles ${\mathcal U}_{y}$
and ${\mathcal U}_{x}$ are non-isomorphic and stable. Hence
$H^0({\mathcal M}_{\xi}, \,{\mathcal U}_{y} \otimes {\mathcal
U}^*_{x})=0.$ This implies that
\begin{equation}\label{4}
p_{X*}({\mathcal U}_{\xi}
\otimes p_{\mathcal M}^*{\mathcal U}^*_{x}) =0\, .
\end{equation}
So \eqref{3} gives
$$
{\mathcal R}^0\,=\,{\mathcal R}^0 p_{3*}(p_{12}^* {\mathcal
U}_{\xi} \otimes p^*_1L_1\otimes p^*_{23}{\mathcal U}^*_\xi) \,=\,
0
$$
and the proof is complete.
\end{proof}

In the following sections we will use Hecke transformations to
compute the fibre of ${\mathcal R}^i$ for $i=1,2$, and will prove
the following propositions.

\begin{propn}\label{2.4}
${\mathcal R}^2 \,=\,0.$
\end{propn}

\begin{propn}\label{2.5}
${\mathcal R}^1$ is a line bundle. Moreover,
${\mathcal R}^1 \cong L_1\otimes TX.$
\end{propn}

Assume for the moment Propositions \ref{2.4} and \ref{2.5}.

{}From the exact sequence \eqref{1} and the previous propositions  we
have the following theorem.

\begin{thm}\label{2.6}
$H^0({\mathcal M}_{\xi},{\mathcal W}_{\xi}(L_1)\otimes
{\mathcal W}_{\xi}(L_2)^*)\cong H^0(X,L_1\otimes L^*_2).$
\end{thm}

\begin{proof}
Combining \eqref{1}, Proposition \ref{2.3} and the second part
of Proposition \ref{2.1} it follows that
$$ H^0({\mathcal M}_{\xi},\, {\mathcal W}_{\xi}(L_1)\otimes
{\mathcal W}_{\xi}(L_2)^*) \, \cong\, H^0(X,\, {\mathcal R}^1\otimes
 L^*_2\otimes K_X)\, .
$$
{}From Proposition \ref{2.5} it follows immediately that
$$
H^0(X,\, {\mathcal R}^1\otimes L^*_2\otimes K_X)
\, \cong \, H^0(X,\, L_1\otimes L^*_2)
$$
and hence the proof is complete.
\end{proof}

\begin{coro}\label{2.7}
If $L_1 \,\not\cong \, L_2$, then ${\mathcal W}_{\xi}(L_1)
\, \not\cong \, {\mathcal W}_{\xi}(L_2)$.
\end{coro}

\begin{proof}
Since $L_1\not\cong L_2$, we have
$$
H^0({\mathcal M}_{\xi},\, {\mathcal W}_{\xi}(L_1)\otimes
{\mathcal W}_{\xi}(L_2)^*)\,=\, 0\, .
$$
Consequently, $ {\mathcal W}_{\xi}(L_1)\not\cong
{\mathcal W}_{\xi}(L_2)$.
\end{proof}

The following is also a corollary of Theorem \ref{2.6}.

\begin{coro}\label{2.8}
The vector bundle ${\mathcal W}_{\xi}(L)$ is simple
for any $L\in J$. In other words,
$H^0({\mathcal M}_{\xi},{\mathcal W}_{\xi}(L)\otimes
{\mathcal W}_{\xi}(L)^*)\cong {\mathbb C}$.
\end{coro}

The following theorem gives the infinitesimal deformations
of ${\mathcal W}_{\xi}(L)$.

\begin{thm}\label{2.9}
For any line bundle $L\,\in\, J$,
the space of infinitesimal deformations of the vector bundle
${\mathcal W}_{\xi}(L)$, namely $H^1({\mathcal M}_\xi ,\, End
( {\mathcal W}_{\xi}(L)))$, is canonically isomorphic to
$H^1(X,\, {\mathcal O}_X)$. In particular,
$$
\dim H^1({\mathcal M}_\xi ,\, End  ( {\mathcal W}_{\xi}(L)))
\, =\, g\, .
$$
\end{thm}

\begin{proof}
Let $L_1=L_2=L$. From \eqref{2} and Proposition \ref{2.4} we
obtain an isomorphism
\begin{equation}\label{5}
H^1(X,\, {\mathcal R}^1\otimes L^* \otimes
 K_X)\, \cong H^2(X\times {\mathcal M}_{\xi}\times X,
p^*_{12}{\mathcal U}_{\xi}\otimes p^*_1L \otimes p^*_{23}
{\mathcal U}_{\xi}^* \otimes p^*_3 L^*\otimes p^*_3 K_X)\, .
\end{equation}

{}From Proposition \ref{2.5} we have
$$\begin{array}{lll}
H^1(X,\, {\mathcal R}^1 \otimes
L^*\otimes K_X)\,& \cong& H^1(X,L\otimes TX \otimes L^*\otimes K_X)\\
&\cong &H^1(X, \, {\mathcal O}_X)\, .
\end{array}
$$
Combining this observation with Proposition \ref{2.1}
and \eqref{5} we get
$$
H^1({\mathcal M}_\xi ,\, End ( {\mathcal W}_{\xi}(L)))
\,\cong \, H^1(X, \, {\mathcal O}_X)\, .
$$
\end{proof}

\begin{rem}\label{2.10}
\begin{em} From the proof of Theorem
\ref{2.9} and Proposition \ref{2.1} we see that if $L_1$ and $L_2$
are not isomorphic then
$$H^1({\mathcal M}_\xi ,\, {\mathcal W}_{\xi}(L_1)
\otimes {\mathcal W}_{\xi}(L_2)^*)
\,\cong \, H^1(X, \, L_1\otimes L^*_2)\, .$$
Hence
$${\rm dim}H^1({\mathcal M}_\xi ,\, {\mathcal W}_{\xi}(L_1)
\otimes {\mathcal W}_{\xi}(L_2)^*)\, =\, g-1\, .$$
\end{em}
\end{rem}

\section{The Hecke transformation}

In this section we will use Hecke transformations to compute
the cohomology groups $H^i(X\times{\mathcal M}_\xi ,\,
{\mathcal U}_{\xi}\otimes p^*_1L_1\otimes
p_{\mathcal M}^*{\mathcal U}^*_x)$ for any $x\in X$ and
prove Proposition \ref{2.4}.
The details of the Hecke transformation and its properties
can be found in \cite{nr,nr2}.
We will briefly describe it and
note those properties that will be needed here.

Fix a point $x\in X$. Let ${\mathbb P}({\mathcal U}_x)$ denote the
projective bundle over ${\mathcal M}_\xi$ consisting of lines in
${\mathcal U}_x$. If $f$ denotes the natural projection of
${\mathbb P}({\mathcal U}_x)$ to ${\mathcal M}_\xi$ and
${\mathcal O}_{{\mathbb P}({\mathcal U}_x)}(-1)$  the
tautological line bundle then
$$f_*{\mathcal O}_{{\mathbb P}({\mathcal U}_x)}(1) \,
\cong\, {\mathcal U}^*_x\, ,
$$
and ${\mathcal R}^j f_* {\mathcal O}_{{\mathbb P}
({\mathcal U}_x)}(1) \, = 0$
for all $j>0.$ From the commutative diagram
$$
\begin{array}{ccc}
X\times {\mathbb P}({\mathcal U}_x)&
\stackrel{{\rm Id}_X\times f}{\longrightarrow}&X \times {\mathcal M}_\xi\\
p_{{\mathbb P}({\mathcal U}_x)}\Big\downarrow&&\Big\downarrow p_{\mathcal
M}\\ {\mathbb P}({\mathcal U}_x)&
\stackrel{f}{\longrightarrow}&{\mathcal M}_\xi
\end{array}
$$
and the base change theorem, we
 deduce that
$$H^i(X\times{\mathcal M}_\xi ,\,
{\mathcal U}_{\xi}\otimes p^*_XL_1\otimes p_{\mathcal
M}^*{\mathcal U}^*_x)$$
\begin{equation}\label{6}
\cong
H^i(X\times {\mathbb P}({\mathcal U}_x),\, ({\rm Id}_X\times f)^*
({\mathcal U}_{\xi}\otimes p^*_XL_1)
\otimes p_{{\mathbb P}({\mathcal U}_x)}^*
{\mathcal O}_{{\mathbb P}({\mathcal
U}_x)}(1))
\end{equation}
for all $i$.

Moreover, since $p_{{\mathbb P}({\mathcal U}_x)*}({\rm Id}_X\times f)^*
({\mathcal U}_{\xi}\otimes p^*_XL_1) \, \cong\, f^*{\mathcal
W}_\xi(L_1)$, there is a canonical isomorphism
$$
H^i(X\times {\mathbb P}({\mathcal U}_x),\, ({\rm Id}_X\times f)^*
({\mathcal U}_{\xi}\otimes p^*_XL_1 )
\otimes p_{{\mathbb P}({\mathcal U}_x)}^*{\mathcal O}_{{\mathbb
P}({\mathcal U}_x)}(1))$$
\begin{equation}\label{7}
\cong\,
H^i({\mathbb P}({\mathcal U}_x),\, f^*{\mathcal
W}_\xi(L_1)\otimes {\mathcal O}_{{\mathbb P}({\mathcal U}_x)}(1))
\end{equation}
for all $i$.

To compute the cohomology groups
$H^i({\mathbb P}({\mathcal U}_x),\, f^*{\mathcal
W}_\xi(L_1)\otimes {\mathcal O}_{{\mathbb P}({\mathcal U}_x)}(1))$
we use Hecke transformations.

A point in ${\mathbb P}({\mathcal U}_x)$
represents a vector bundle
$E$ and a line $l$ in the fiber $E_x$ at $x$, or equivalently a
non-trivial exact sequence
\begin{equation}\label{8}
0\, \longrightarrow\, E\, \longrightarrow\, F \,
\longrightarrow\, {\mathbb C}_x \, \longrightarrow\, 0 \,
\end{equation}
determined up to a scalar multiple; here ${\mathbb C}_x$ denotes the
torsion sheaf supported at $x$ with stalk $\mathbb C$. The sequences
\eqref{8} fit together to form a universal sequence
\begin{equation}\label{9}
0\, \longrightarrow\, ({\rm Id}_X\times f)^*
{\mathcal U}_{\xi}\otimes p_{{\mathbb P}({\mathcal U}_x)}^*{\mathcal
O}_{{\mathbb P}({\mathcal U}_x)}(1)\, \longrightarrow\, {\mathcal F} \,
\longrightarrow\, p_X^*{\mathbb C}_x \, \longrightarrow\, 0 \,
\end{equation}
on $X\times{\mathbb P}({\mathcal U}_x)$.
  If $\eta$ denotes the line bundle $\xi\otimes {\mathcal O}_X(
x)$ over $X$ and ${\mathcal
M}_\eta $ the moduli space of stable bundles $ {\mathcal M}_\eta(n,d+1)$
then from \eqref{8} and \eqref{9} we get a rational map
$$
\gamma\, :\, {\mathbb P}({\mathcal U}_x)\, - - \rightarrow\,{\mathcal
M}_\eta$$
which sends any pair $(E,l)$ to $F$. This map is not everywhere
defined since the bundle $F$ in \eqref{8} need not be stable.

Our next
object is to find a Zariski-open subset
$Z$ of ${\mathcal M}_\eta$, over which $\gamma$ is defined and is a
projective fibration, such that the complement of $Z$ in ${\mathcal
M}_\eta$ has codimension at least 4. The construction and
calculations are similar to those of \cite[Proposition 6.8]{nr}, but
our results do not seem to follow directly from that proposition.

As in \cite[\S8]{nr} or \cite[\S5]{nr2}, we define a bundle $F$ to
be $(0,1)$-{\it stable} if, for every proper subbundle $G$ of $F$,
$$\frac{\deg G}{\hbox{rk}G}<\frac{\deg F-1}{\hbox{rk}F}.$$
Clearly every $(0,1)$-stable bundle is stable. We denote by $Z$ the
subset of ${\mathcal M}_{\eta}$ consisting of $(0,1)$-stable bundles.

\begin{lemma}\label{3.1}
{\rm (i)} $Z$ is a Zariski-open subset of ${\mathcal
M}_{\eta}$ whose complement has codimension at least $4$.

{\rm (ii)} $\gamma$ is a projective fibration over $Z$ and
$\gamma^{-1}(Z)$ is a Zariski-open subset of ${\mathbb P}({\mathcal
U}_x)$ whose complement has codimension at least $4$.
\end{lemma}

\begin{proof}
(i) The fact that $Z$ is Zariski-open is standard (see
\cite[Proposition 5.3]{nr2}).

The bundle $F\in{\mathcal M}_{\eta}$ of rank $n$ and degree $d+1$
fails to be $(0,1)$-stable if and only if it has a subbundle $G$ of
rank $r$ and degree $e$ such that $ne\ge r((d+1)-1)$, i.e.,
\begin{equation}\label{10}
rd\,\le\, ne\, .
\end{equation}
By considering the extensions
$$0\longrightarrow G\longrightarrow F\longrightarrow H\longrightarrow
0,$$
we can estimate the codimension of ${\mathcal M}_{\eta}-Z$ and show
that it is at least
\begin{equation}\label{11}
\delta=r(n-r)(g-1)+(ne-r(d+1))
\end{equation}
(compare the proof of \cite[Proposition 5.4]{nr2}). Note that, since
$(n,d)=1$, \eqref{10} implies that $rd\le ne-1$. Given that $g\ge3$, we
see that $\delta<4$ only if $g=3$, $n=2,3$ or $g=4$, $n=2$. These are
exactly the cases that were excluded in the introduction.

(ii) $\gamma^{-1}(Z)$ consists of all pairs $(E,l)$ for which the
bundle $F$ in \eqref{8} is $(0,1)$-stable. As in (i), this is a
Zariski-open subset.
It follows at once from \eqref{10} that, if $F$ is $(0,1)$-stable,
then $E$ is stable. So, if $F\in Z$, it follows from
\eqref{8} that $\gamma^{-1}(F)$ can be identified with the projective
space ${\mathbb P}(F_x^*)$. Using the universal projective bundle on
$X\times{\mathcal M}_{\eta}$, we see that $\gamma^{-1}(Z)$ is a
projective fibration over $Z$ (not necessarily locally trivial).

Suppose now that $(E,l)$ belongs to the complement of
$\gamma^{-1}(Z)$ in ${\mathbb P}({\mathcal U}_x)$. This means that the
bundle $F$ in \eqref{8} is not $(0,1)$-stable and therefore possesses a
subbundle $G$ satisfying \eqref{10}. If $G\subset E$, this contradicts the
stability of $E$. So there exists an exact sequence
$$0\longrightarrow G'\longrightarrow G\longrightarrow{\mathbb
C}_x\longrightarrow 0$$
with $G'$ a subbundle of $E$ of rank $r$ and degree $e-1$. Moreover,
since $G$ is a subbundle of $F$, $G'_x$ must contain the line $l$.
For fixed $r,e$, these conditions determine a subvariety of ${\mathbb
P}({\mathcal U}_x)$ of dimension at most
$$
(r^2(g-1)+1)+((n-r)^2(g-1)+1)-g+(r-1)
+((g-1)r(n-r)+(rd-n(e-1))-1).
$$
Since $\dim{\mathbb
P}({\mathcal U}_x)=n^2(g-1)-g+n$, a simple calculation shows that the
codimension is at least the number $\delta$ given by \eqref{11}. As in
(i),
this gives the required result.
\end{proof}

By Lemma \ref{3.1}(ii) and a Hartogs-type theorem (see \cite[Theorem
3.8 and Proposition 1.11]{hart}) we have an isomorphism
\begin{equation}\label{12}
H^i({\mathbb P}({\mathcal U}_x),\, f^*{\mathcal
W}_\xi(L_1)\otimes {\mathcal O}_{{\mathbb P}({\mathcal U}_x)}(1))
\cong H^i(\gamma ^{-1}(Z),\, f^*{\mathcal
W}_\xi(L_1)\otimes {\mathcal O}_{{\mathbb P}
({\mathcal U}_x)}(1)|_{\gamma^{-1}(Z)})
\end{equation}
for $i\le2$.

Now let $F\in Z$. As in the proof of Lemma \ref{3.1}, we identify
$\gamma^{-1}(F)$ with ${\mathbb P}(F_x^*)$ and denote it by ${\mathbb P}$.
On $X\times{\mathbb P}$ there is a universal exact sequence
\begin{equation}\label{13}
0\longrightarrow{\mathcal E}\longrightarrow
p_X^*F\longrightarrow p_{\mathbb P}^*{\mathcal O}_{\mathbb P}(1)\otimes
p_X^*{\mathbb C}_x\longrightarrow 0\, .
\end{equation}
The
restriction of \eqref{13} to any point of ${\mathbb P}$ is isomorphic to
the
corresponding sequence \eqref{8}.

\begin{propn}\label{3.2}
Let ${\mathcal F}$ be defined by the universal sequence
\eqref{9}. Then
$${\mathcal F}\big\vert_{X\times{\mathbb P}}\,\cong\, p_X^*F\otimes
p_{\mathbb
P}^*{\mathcal O}_{\mathbb P}(-1)\, .$$
\end{propn}

\begin{proof}
Restricting \eqref{9} to $X\times{\mathbb P}$
gives
$$
0\, \longrightarrow\, ({\rm Id}_X\times f)^*
{\mathcal U}_{\xi}\otimes p_{{\mathbb P}({\mathcal U}_x)}^*{\mathcal
O}_{{\mathbb P}({\mathcal U}_x)}(1))\big\vert_{X\times{\mathbb P}}
\longrightarrow\, {\mathcal F}\big\vert_{X\times{\mathbb P}}
\,
\longrightarrow\, p_X^*{\mathbb C}_x \, \longrightarrow\, 0.
$$
This must coincide with the universal sequence
\eqref{13} up to tensoring
by some line bundle lifted from ${\mathbb P}$. The result follows.
\end{proof}

Next we tensor \eqref{9} by $p_X^*L_1$, restrict it to
$X\times\gamma^{-1}(Z)$ and take the direct image on
$\gamma^{-1}(Z)$. This gives
\begin{equation}\label{14}
0\, \longrightarrow\, f^*
{\mathcal W}_{\xi}(L_1)\otimes {\mathcal O}_{{\mathbb P}({\mathcal
U}_x)}(1)\big\vert_{\gamma^{-1}(Z)} \longrightarrow\, p_{{\mathbb
P}({\mathcal U}_x)*}({\mathcal F}\otimes
p_X^*L_1)\big\vert_{\gamma^{-1}(Z)}
\,
\longrightarrow\, {\mathcal O}_{\gamma^{-1}(Z)} \, \longrightarrow\,
0\, .
\end{equation}

\begin{propn}\label{3.3}
$R^i_{\gamma*}(p_{{\mathbb P}({\mathcal
U}_x)*}({\mathcal F}\otimes p_X^*L_1)\big\vert_{\gamma^{-1}(Z)})=0$
for all $i$.
\end{propn}

\begin{proof}
It is sufficient to show that $p_{{\mathbb P}({\mathcal
U}_x)*}({\mathcal F}\otimes p_X^*L_1)\big\vert_{\mathbb P}$ has trivial
cohomology. By Proposition \ref{3.2}
\begin{eqnarray*}p_{{\mathbb P}({\mathcal U}_x)*}({\mathcal F}\otimes
p_X^*L_1)\big\vert_{\mathbb P}&\cong& p_{{\mathbb
P}*}(p_X^*F\otimes p_X^*L_1\otimes p_{\mathbb P}^*{\mathcal
O}_{\mathbb P}(-1))\\&\cong& H^0(X,F\otimes
L_1)\otimes{\mathcal O}_{\mathbb P}(-1)\end{eqnarray*}
and the result follows.
\end{proof}

\begin{coro}\label{3.4}
${\mathcal R}^i_{\gamma_*} (f^*{\mathcal
W}_\xi(L_1)\otimes {\mathcal O}_{{\mathbb P}({\mathcal
U}_x)}(1)\big\vert_{\gamma^{-1}(Z)})
\, = 0$ for $i \ne1$. Moreover,
$${\mathcal R}^1_{\gamma_*} (f^*{\mathcal
W}_\xi(L_1)\otimes {\mathcal O}_{{\mathbb P}({\mathcal
U}_x)}(1)\big\vert_{\gamma^{-1}(Z)})
\cong{\mathcal O}_Z.$$
\end{coro}

\begin{proof}
This follows at once from \eqref{14} and Proposition \ref{3.3}.
\end{proof}

Now we are in a position to compute the cohomology groups of
$H^i(X\times{\mathcal M}_\xi ,\,
{\mathcal U}_{\xi}\otimes p^*_XL_1\otimes
p_{\mathcal M}^*{\mathcal U}^*_x)$ for  $i=1,2.$

\begin{propn}\label{3.5}
$H^2(X\times{\mathcal M}_\xi ,\,
{\mathcal U}_{\xi}\otimes p^*_XL_1\otimes
p_{\mathcal M}^*{\mathcal U}^*_x)=0 $ for any $x\in X.$
\end{propn}

\begin{proof}
The combination of \eqref{6}, \eqref{7} and \eqref{12}
yields
$$
H^2(X\times{\mathcal M}_\xi ,\,
{\mathcal U}_{\xi}\otimes p^*_XL_1\otimes
p_{\mathcal M}^*{\mathcal U}^*_x) \,\cong\,
H^2(\gamma ^{-1}(Z),\, f^*{\mathcal
W}_\xi(L_1)\otimes {\mathcal O}_{{\mathbb P}
({\mathcal U}_x)}(1)\big\vert_{\gamma^{-1}(Z)})\, .
$$
Using Corollary \ref{3.4} and Lemma \ref{3.1}(i), the Leray spectral
sequence
for the map $\gamma$ gives$$
\begin{array}{lll}
H^2(\gamma^{-1}(Z),\,   f^*{\mathcal
W}_\xi(L_1)\otimes {\mathcal O}_{{\mathbb P}
({\mathcal U}_x)}(1)\big\vert _{\gamma^{-1}(Z)})\,
&\cong\,& H^{1}(Z,\, {\mathcal O}_Z)\\
&\cong\,& H^{1}({\mathcal M}_\eta,\, {\mathcal O}_{{\mathcal
M}_{\eta}})
\end{array}
$$
It is known that
$H^1({\mathcal M}_\eta,\, {\mathcal O}_{{\mathcal M}_{\eta}})\,=\,0$
\cite{dn}. Therefore,
$$
H^2(X\times{\mathcal M}_\xi ,\,
{\mathcal U}_{\xi}\otimes p^*_XL_1\otimes
p_{\mathcal M}^*{\mathcal U}^*_x) = 0\, .
$$
\end{proof}

{\it Proof of Proposition \ref{2.4}.}\, Proposition \ref{2.4} is
an immediate consequence of Proposition \ref{3.5}.
$\hfill{\Box}$
\medskip

\begin{propn}\label{3.6}
For any point $x\in X$, $\dim H^1(X\times{\mathcal M}_\xi ,\,
{\mathcal U}_{\xi}\otimes p^*_XL_1\otimes
p_{\mathcal M}^*{\mathcal U}^*_x)\, =\, 1$.
\end{propn}

\begin{proof}
As in the proof of Proposition \ref{3.5}
we conclude that
$$
H^1(X\times{\mathcal M}_\xi ,\,
{\mathcal U}_{\xi}\otimes p^*_XL_1\otimes
p_{\mathcal M}^*{\mathcal U}^*_x)
\, \cong\,H^{0}({\mathcal M}_\eta,\, {\mathcal O}_{{\mathcal
M}_{\eta}})\, .
$$
Now ${\mathcal M}_{\eta}$ is just the non-singular part of the moduli
space of semistable bundles of rank $n$ and determinant $\eta$, and
the latter space is complete and normal. So
$\dim H^0({\mathcal M}_{\eta},\, {\mathcal O}_{{\mathcal
M}_{\eta}})\,=\, 1$.
\end{proof}

\begin{rem}\label{3.7}
\begin{em}
Since the fibres of $\gamma$ are projective spaces, we have
$\gamma_*{\mathcal O}_{\gamma^{-1}(Z)}\cong{\mathcal O}_Z$ and all the
higher direct images of ${\mathcal O}_{\gamma^{-1}(Z)}$ are $0$. Hence
$$H^i(Z,{\mathcal O}_Z)\cong H^i(\gamma^{-1}(Z),{\mathcal
O}_{\gamma^{-1}(Z)})$$for all $i$. Similarly
$$H^i({\mathbb P}({\mathcal U}_x),{\mathcal O}_{{\mathbb P}({\mathcal
U}_x)})\cong H^i({\mathcal M}_{\xi},{\mathcal O}_{{\mathcal
M}_{\xi}})=0$$ for
$i>0$ since
${\mathcal M}_{\xi}$ is a smooth projective rational variety. It
follows from the proof of Lemma \ref{3.1} that,
if we define $\delta$ as in
\eqref{11},
$$\delta\ge i+2\ge3\Rightarrow H^i(Z,{\mathcal O}_Z)=0.$$
The proof of Proposition \ref{3.5} now gives
$$\delta\ge i+2\ge4\Rightarrow H^i(X\times{\mathcal M}_{\xi},
{\mathcal U}_{\xi}\otimes p^*_XL_1\otimes p^*_{\mathcal M}{\mathcal
U}_x^*)=0.$$
This in turn implies that $R^i=0$. Proposition \ref{2.1} and the Leray
spectral sequence of $p_3$ (cf. \eqref{1} and \eqref{2}) now give
$$\delta\ge i+3\ge5\Rightarrow H^i({\mathcal M}_{\xi},{\mathcal
W}_{\xi}(L_1)\otimes {\mathcal W}_{\xi}(L_2)^*)=0.$$
In particular, taking $i=2$ and $L_1=L_2=L$, and using
\eqref{11}, we obtain
$$H^2({\mathcal M}_{\xi},End({\mathcal W}_{\xi}(L)))=0$$
except possibly when $g=3, n=2, 3, 4$; $g=4,n=2$; $g=5,n=2$.
\end{em}
\end{rem}

\section{Proof of Proposition \ref{2.5}}

Let $\Delta $ be the diagonal divisor in $X\times X$. Pull back the
exact sequence
$$0\, \longrightarrow\, {\mathcal O}(-{\Delta}) \, \longrightarrow\,
{\mathcal O}\, \longrightarrow\, {\mathcal O}_{\Delta}\,
\longrightarrow\, 0
$$
to $X\times {\mathcal M}_{\xi} \times X$ and tensor it with
$ p^*_{12}
{\mathcal U}_{\xi}\otimes p^*_1L_1 \otimes
p^*_{23}{\mathcal U}^*_\xi$. Now, the direct image
sequence for the projection
$p_3$ gives the following exact sequence over $X$
$$
\longrightarrow\, {\mathcal R}^i p_{3*} ( p^*_{12}
{\mathcal U}_{\xi}\otimes p^*_1L_1 \otimes
p^*_{23}{\mathcal U}^*_\xi (-{\Delta}))\, \longrightarrow\,
{\mathcal R}^i
 \, \longrightarrow\, {\mathcal R}^i p_{3*} ( p^*_{12}
{\mathcal U}_{\xi}\otimes p^*_1L_1 \otimes
p^*_{23}{\mathcal U}^*_\xi |_{\Delta \times {\mathcal M}_{\xi}})
$$
\begin{equation}\label{15}
\longrightarrow\, {\mathcal R}^{i+1} p_{3*} ( p^*_{12}
{\mathcal U}_{\xi}\otimes p^*_1L_1 \otimes
p^*_{23}{\mathcal U}^*_\xi (-{\Delta}))
\,\longrightarrow\,\ldots
\end{equation}

The following propositions will be used in computing the direct
images ${\mathcal R}^i.$

\begin{propn}\label{4.1}
For any $L_1\in J$, the direct images of
$$
p^*_{12}
{\mathcal U}_{\xi}\otimes p^*_1L_1 \otimes
p^*_{23}{\mathcal U}^*_\xi |_{\Delta \times {\mathcal M}_{\xi}}
$$
have the following description:

\begin{enumerate}
\item ${\mathcal R}^0 p_{3*}( p^*_{12}
{\mathcal U}_{\xi}\otimes p^*_1L_1 \otimes
p^*_{23}{\mathcal U}^*_\xi |_{\Delta \times {\mathcal M}_{\xi}})
\,\cong\, L_1$
\item ${\mathcal R}^1 p_{3*}( p^*_{12}
{\mathcal U}_{\xi}\otimes p^*_1L_1 \otimes
p^*_{23}{\mathcal U}^*_\xi |_{\Delta \times {\mathcal M}_{\xi}})
\,\cong\, L_1\otimes TX$
\item  ${\mathcal R}^2 p_{3*}( p^*_{12}
{\mathcal U}_{\xi}\otimes p^*_1L_1 \otimes
p^*_{23}{\mathcal U}^*_\xi |_{\Delta \times {\mathcal M}_{\xi}})\,=\,0$
\end{enumerate}
where $TX$ is the tangent bundle of $X$.
\end{propn}

\begin{proof}
Identifying $\Delta $ with $X$ we have
$${\mathcal R}^i p_{3*} ( p^*_{12}
{\mathcal U}_{\xi}\otimes p^*_1L_1 \otimes
p^*_{23}{\mathcal U}^*_\xi |_{\Delta \times {\mathcal M}_{\xi}})
\,\cong \, {\mathcal R}^i p_{X*}(
{\mathcal U}_{\xi}\otimes p^*_XL_1 \otimes
{\mathcal U}^*_\xi ).$$

The proposition follows from a result of
Narasimhan and Ramanan \cite[Theorem 2]{nr} that says
\begin{equation}\label{16}
H^i({\mathcal M}_{\xi},\, {\mathcal U}_x \otimes {\mathcal U}^*_x )\,
\cong \,
\left\{\begin{array}{ll} {\mathbb C} &\mbox{if} \ \ i=0,1\\
0 & \mbox{if} \ \
i=2.  \end{array}\right.
\end{equation}
For $i=0$ the isomorphism is given
by the obvious inclusion of ${\mathcal O}_{{\mathcal M}_{\xi}}$ in
${\mathcal U}_x
\otimes {\mathcal U}^*_x $ and therefore globalises to give
${\mathcal R}^0 p_{X*}(  {\mathcal U}_{\xi} \otimes
{\mathcal U}^*_\xi )\cong{\mathcal O}_X$. Similarly for $i=1$ the
isomorphism is given by the infinitesimal deformation map of
${\mathcal U}_{\xi}$ regarded as a family of bundles over $\mathcal
M_{\xi}$ parametrised by $X$; this globalises to ${\mathcal R}^1
p_{X*}({\mathcal U}_{\xi}\otimes
{\mathcal U}^*_\xi )\cong TX$.
\end{proof}

Propositions \ref{4.1}, \ref{2.3} and the exact sequence \eqref{15}
together give the following exact sequence of direct images
$$
0\,\longrightarrow\, L_1 \,\longrightarrow\, {\mathcal R}^1 p_{3*}(
p^*_{12}{\mathcal U}_{\xi}\otimes p^*_1L_1 \otimes
p^*_{23}{\mathcal U}^*_\xi (-{\Delta}))\, \longrightarrow\,
{\mathcal R}^1 \, \stackrel{\alpha}{\longrightarrow}\,
L_1\otimes TX
$$
\begin{equation}\label{17}
\longrightarrow \,  {\mathcal R}^2 p_{3*}(
p^*_{12}{\mathcal U}_{\xi}\otimes p^*_1L_1 \otimes
p^*_{23}{\mathcal U}^*_\xi (-{\Delta}))
\longrightarrow\ldots
\end{equation}
For any $x\in X$ we have the cohomology exact sequence
$$
\longrightarrow\, H^i(X\times {\mathcal M}_{\xi}, \,{\mathcal U}_{\xi}
\otimes p^*_XL_1 \otimes
p^*_{\mathcal M}{\mathcal U}^*_x(-x)) \, \longrightarrow\,
 H^i(X\times {\mathcal M}_{\xi}, \, {\mathcal U}_{\xi}
\otimes p^*_XL_1 \otimes
p^*_{\mathcal M}{\mathcal U}^*_x)
$$
\begin{equation}\label{18}
\longrightarrow\, H^i({\mathcal M}_{\xi} ,\,
{\mathcal U}_x\otimes {\mathcal U}^*_x)\otimes (L_1)_x
\, \longrightarrow\,
H^{i+1}(X\times {\mathcal M}_{\xi}, \,{\mathcal U}_{\xi}
\otimes p^*_XL_1 \otimes p^*_{\mathcal M}{\mathcal U}^*_x(-x))\, ,
\end{equation}
where $(L_1)_x$ is the fiber of $L_1$ at $x$.

By \eqref{4}, $p_{X*} ({\mathcal U}_{\xi}\otimes
p_{\mathcal M}^*{\mathcal U}^*_x) =0$. So the Leray spectral
sequence for $p_X$ gives
$$
H^1(X\times{\mathcal M}_\xi ,\,
{\mathcal U}_{\xi}\otimes p^*_XL_1\otimes
p_{\mathcal M}^*{\mathcal U}^*_x)\,\cong\, H^0(X,\, {\mathcal
R}^1 p_{X*} ({\mathcal U}_{\xi}\otimes p^*_XL_1\otimes
p_{\mathcal M}^*{\mathcal U}^*_x))
$$
and
$$
H^1(X\times{\mathcal M}_\xi ,\,
{\mathcal U}_{\xi}\otimes p^*_XL_1\otimes
p_{\mathcal M}^*{\mathcal U}^*_x(-x))\,\cong\, H^0(X,\,
{\mathcal R}^1 p_{1*} ({\mathcal U}_{\xi}\otimes p^*_XL_1\otimes
p_{\mathcal M}^*{\mathcal U}^*_x)(-x)).
$$
Since ${\mathcal U}_x$ is simple \cite[Theorem 2]{nr},
\eqref{18} gives the exact sequence
$$
0\,\longrightarrow\, (L_1)_x\, \longrightarrow \, H^0(X,{\mathcal
R}^1 p_{X*}({\mathcal U}_{\xi}\otimes p^*_XL_1\otimes
p_{\mathcal M}^*{\mathcal U}^*_x)(-x))\hspace{3.cm}
$$
$$
\hspace{3.cm} \longrightarrow\,
H^0(X,{\mathcal R}^1 p_{X*}
({\mathcal U}_{\xi}\otimes p^*_XL_1\otimes
p_{\mathcal M}^*{\mathcal U}^*_x)) \, \longrightarrow\ldots
$$
This implies that ${\mathcal R}^1 p_{X*}
({\mathcal U}_{\xi}\otimes p^*_XL_1\otimes
p_{\mathcal M}^*{\mathcal U}^*_x)$ has torsion at $x$. Now from
\eqref{16} we conclude that ${\mathcal R}^1 p_{X*}
({\mathcal U}_{\xi}\otimes p^*_XL_1\otimes
p_{\mathcal M}^*{\mathcal U}^*_x)$
is a torsion sheaf, and hence
$$
H^1(X,\, {\mathcal R}^1 p_{X*}
({\mathcal U}_{\xi}\otimes p^*_XL_1\otimes
p_{\mathcal M}^*{\mathcal U}^*_x)(-x))\, =\, H^1(X,\,
{\mathcal R}^1 p_{X*} ({\mathcal U}_{\xi}\otimes p^*_XL_1\otimes
p_{\mathcal M}^*{\mathcal U}^*_x))\,=\,0\, .
$$

The
Leray spectral sequence for $p_X$ now yields
\begin{equation}\label{19}
H^2(X\times{\mathcal M}_\xi ,\,
{\mathcal U}_{\xi}\otimes p^*_XL_1\otimes
p_{\mathcal M}^*{\mathcal U}^*_x)\, \cong \, H^0(X,\, {\mathcal
R}^2 p_{X*} ({\mathcal U}_{\xi}\otimes p^*_XL_1\otimes
p_{\mathcal M}^*{\mathcal U}^*_x))
\end{equation}
and
\begin{equation}\label{20}
H^2(X\times{\mathcal M}_\xi ,\,
{\mathcal U}_{\xi}\otimes p^*_XL_1\otimes
p_{\mathcal M}^*{\mathcal U}^*_x(-x))\cong H^0(X,{\mathcal R}^2
p_{X*} ({\mathcal U}_{\xi}\otimes p^*_XL_1\otimes
p_{\mathcal M}^*{\mathcal U}^*_x)(-x))\, .
\end{equation}
Now from \eqref{16} it follows that
${\mathcal R}^2 p_{X*}
({\mathcal U}_{\xi}\otimes p^*_XL_1\otimes
p_{\mathcal M}^*{\mathcal U}^*_x)$ is a torsion sheaf, and from
Proposition \ref{3.5} and \eqref{19} that its space of sections is $0$.
So ${\mathcal R}^2 p_{X*}
({\mathcal U}_{\xi}\otimes p^*_XL_1\otimes
p_{\mathcal M}^*{\mathcal U}^*_x)=0$ and by \eqref{20}
$$
H^2(X\times{\mathcal M}_\xi ,\,
{\mathcal U}_{\xi}\otimes p^*_XL_1\otimes
p_{\mathcal M}^*{\mathcal U}^*_x(-x))\,=\, 0\, .$$

So, we have the following proposition:

\begin{propn}\label{4.2}
$ {\mathcal R}^2 p_{3*}( p^*_{12}
{\mathcal U}_{\xi}\otimes p^*_1L_1 \otimes
p^*_{23}{\mathcal U}^*_\xi (-{\Delta}))\,= \,0.$
\end{propn}

Now Proposition \ref{2.5} is easy to derive.

\medskip
{\it Proof of Proposition \ref{2.5}.}\, By Proposition \ref{3.6}
$$
H^1(X\times{\mathcal M}_\xi ,\,
{\mathcal U}_{\xi}\otimes p^*_XL_1\otimes
p_{\mathcal M}^*{\mathcal U}^*_x)\cong {\mathbb C}\, .
$$
Therefore, ${\mathcal R}^1 = {\mathcal R}^1 p_{3*}( p^*_{12}
{\mathcal U}_{\xi}\otimes p^*_1L_1 \otimes
p^*_{23}{\mathcal U}^*_\xi)$ is a line bundle.
Moreover Proposition \ref{4.2} implies that  the map
$\alpha \, :\, {\mathcal R}^1 \, \longrightarrow\,
L_1\otimes TX$ in the exact sequence \eqref{17}
is surjective. Therefore $\alpha$ must be an isomorphism,
and the proof is complete.
$\hfill{\Box}$

Earlier Theorems \ref{2.6} and \ref{2.9} were proved assuming
Propositions \ref{2.4} and \ref{2.5}. Therefore Theorems \ref{2.6} and
\ref{2.9} are now established.

\medskip

\section{Family of Deformations}

Fix a point $x\in X$. Let $\mathcal L$ be the Poincar\'e line
bundle over $X\times J$ which is trivial on
$x\times  J$.
It was proved in \cite{bbn} that the family
$\widetilde{{\mathcal U}_{\xi}}:=
p^*_{12}{\mathcal L}\otimes p^*_{13}{\mathcal U}_{\xi}$ over
$X\times J\times {\mathcal M}_{\xi}$ is a complete
family of deformations of ${\mathcal U}_{\xi}$ parametrised by
$J$.
In this section we will prove that the direct image
$$\widetilde{{\mathcal W}}:= p_{23*}(\widetilde{{\mathcal U}_{\xi}})$$
is a complete family of deformations of ${\mathcal W}_{\xi}(L)$ for any
$L\in J$.

First note that $\widetilde{{\mathcal W}}$ is locally free. Indeed,
${\mathcal R}^i p_{23*}(\widetilde{{\mathcal U}_{\xi}}) =0$ for
$i\not= 0.$ Moreover, for each $L\in J$,
$$\widetilde{{\mathcal W}}\big\vert _{\{L\}\times {\mathcal M}_{\xi}}
\cong  {\mathcal W}_{\xi}(L)\, . $$

In \cite{bbgn} it was proved that ${\mathcal W}_{\xi}$ is stable with
respect to the unique polarisation of ${\mathcal M}_{\xi}$. From the
proof of the Theorem in
\cite{bbgn} we can see that  for any line bundle $L\in J$,
${\mathcal W}_{\xi}(L)$ is also stable.
Since  $J$ is
irreducible, the Hilbert polynomial
$P$ of ${\mathcal W}_{\xi}(L) $ is well defined and independent of
the choice of $L$.

Denote by ${\mathcal M}^0_X({\mathcal W}_{\xi})$
the connected component of the moduli space of stable
sheaves, with the same Hilbert polynomial as
${\mathcal W}_{\xi}({\mathcal O})$, on ${\mathcal M}_{\xi}$
containing  ${\mathcal W}_{\xi}({\mathcal O})$. As in \cite{ib} the
determinant line bundle $M$ on ${\mathcal M}^0_X({\mathcal W}_{\xi})$
defines a polarisation on ${\mathcal M}^0_X({\mathcal W}_{\xi})$.

Define the morphism
$$\phi : J\longrightarrow {\mathcal M}^0_X({\mathcal W}_{\xi})$$
by $L\mapsto {\mathcal W}_{\xi}(L).$

\begin{thm}\label{5.1}
The morphism
$$
\phi\, : \,J\, \longrightarrow\,
{\mathcal M}^0_X({\mathcal W}_{\xi})
$$
is an isomorphism of polarised varieties.
\end{thm}

\begin{proof}
By Corollary \ref{2.7}
$\phi$ is
injective, so its image has dimension $g$. On the other hand, by
Theorem \ref{2.9}, the Zariski tangent space of  ${\mathcal
M}^0_X({\mathcal
W}_{\xi})$ also has dimension $g$ at every point of the image of
$\phi$. It follows that ${\mathcal M}^0_X({\mathcal
W}_{\xi})$ is smooth of dimension $g$ at every point of the image of
$\phi$. Hence, by Zariski's Main Theorem, $\phi$ is an isomorphism
onto an open subset of ${\mathcal M}^0_X({\mathcal
W}_{\xi})$. Finally $J$ is complete and ${\mathcal M}^0_X({\mathcal
W}_{\xi})$ is connected and separated (since ${\mathcal
W}_{\xi}({\mathcal O})$ is stable), so $\phi$ is an isomorphism.

Let $\zeta$ be the polarisation on
${\mathcal M}^0_X({\mathcal W}_{\xi})$ given by the
determinant line bundle \cite[Section 4]{ib}. Let $\Theta$ denote
the principal
polarisation on $J$ defined by a theta divisor.
We wish to show that the isomorphism $\phi$ takes
$\zeta$ to a nonzero constant scalar multiple (independent of the
curve $X$) of $\Theta$.

Take any family of pairs $(X,\xi)$, where $X$ is a connected
non-singular projective curve of genus
$g$ and $\xi$ is a line bundle on $X$ of degree $d>n(2g-2)$,
parametrized  by a connected space
$T$. Consider the corresponding family of moduli spaces ${\mathcal
M}^0_X({\mathcal W}_{\xi})$ (respectively, Jacobians $J$) over $T$,
where
$(X,\xi)$ runs over the family. Using the map $\phi$ an isomorphism
between these two families is obtained.
The polarisation $\zeta$ (respectively, $\Theta$) defines a constant
section of the second direct image over $T$
of the constant sheaf $\mathbb Z$ over the family.
It is known that for the general curve $X$ of genus $g$, the
Neron-Severi group of $J$ is $\mathbb Z$. Therefore,
for such a curve, $\phi$ takes $\zeta$ to a nonzero constant scalar
multiple of $\Theta$. Since $T$ is connected,
if $T$ contains a curve with
$NS(J)\, =\, {\mathbb Z}$, then $\phi$
takes $\zeta$ to the same nonzero constant scalar multiple
of $\Theta$ for every curve in the family. Since the moduli
space of smooth curves of genus $g$ is connected, the proof
is complete.
\end{proof}

\begin{coro}\label{5.2}
The family $\widetilde{{\mathcal W}}$
parametrised by $ J$ is complete. Moreover,
the infinitesimal deformation map of this family at any point of $J$
is an isomorphism.
\end{coro}

\begin{proof}
This follows at once from the theorem.
\end{proof}

\medskip
\begin{rem}\begin {em}
In the proof of Theorem \ref{5.1}, smoothness
follows from the fact that the dimension of ${\mathcal
M}^0_X({\mathcal W}_{\xi})$  is equal to the dimension of its
Zariski tangent space. So we do not need to know that
$H^2({\mathcal M}_{\xi},End({\mathcal W}_{\xi}(L)))=0$
(see Remark  \ref{3.7}) or even that ${\mathcal M}^0_X({\mathcal W}_{\xi})$ is
reduced.
\end{em}\end{rem}
\medskip

Finally we have our Torelli theorem.

\begin{coro}\label{5.4}
Let $X$ and $X'$ be two non-singular algebraic curves of
genus $g\geq 3$ and let $\xi$ (respectively $\xi'$) be a line bundle
of degree $d>n(2g-2)$ on $X$ (respectively $X'$). If
${\mathcal M}^0_X({\mathcal W}_{\xi})\cong
 {\mathcal M}^0_{X'}({\mathcal W}_{\xi'})$ as polarised varieties
then
$X\cong X'$.
\end{coro}

\begin{proof}
This follows at once from Theorem \ref{5.1} and the
classical Torelli theorem.
\end{proof}

%%%%%%%%%%%%%%%%%%%%%%%%%%%%%%%%%%%%%%%%%%%%%%%%%%%%%%%%%%%%%%%%

\end{document}